\documentclass[12pt]{article}
\usepackage{epsfig}
\usepackage{amssymb}

\newtheorem{theorem}{Theorem}[section]
\newtheorem{lemma}[theorem]{Lemma}
\newtheorem{proposition}[theorem]{Proposition}

\newtheorem{definition}[theorem]{Definition}

\newtheorem{conjecture}[theorem]{Conjecture}
\newtheorem{question}[theorem]{Question}

\renewcommand{\H}{\mathbf{H}}

\newenvironment{proof}{{\it Proof:\/}}{$\Box$\vskip 0.08in}
\newcommand{\mod}{{\mbox{ mod }}}

\begin{document}

\begin{center}
\bigskip

{\LARGE \baselineskip=10pt {\ Rational moves and tangle embeddings:
$(2,2)$-moves as a case study} } %
\centerline{ Tokyo, January 29, 2005}

Mieczys{\l }aw K.~D{\c a}bkowski, Makiko Ishiwata and J\'ozef H.~Przytycki
\end{center}

\vspace{27pt} Abstract.\newline
{\footnotesize {\ 
We classify $3$-braids up to $(2,2)$-move equivalence and, in particular, 
we show how to adjust the Harikae-Nakanishi-Uchida conjecture 
so it holds for closed 3-braids.
As an important step in classification of $3$-braids up to 
$(2,2)$-move equivalence
we prove the conjecture for 
$2$-algebraic links and
classify $(2,2)$-equivalence classes for links up to nine crossings.
We also analyze the effect of $(2,2)$-move on Kei (involutive quandle) 
associated to a link.
We construct Burnside Kei, $Q(m,n)$, and ask the question,
motivated by classical Burnside question: for which values
of $m$ and $n$, is $Q(m,n)$ finite?


\ \newline

\section{Introduction\label{1}}

A tangle move is a local modification of a link in which a tangle $T_{1}$ is 
replaced by a tangle $T_{2}$.
The simplest such a move, that reduces every
link in $S^{3}$ into a trivial link, is a crossing change. 
%

It was believed that there are some 
nontrivial \textquotedblleft tangle moves" with the unknotting property,
other than the crossing change. 
For example, the Montesinos-Nakanishi conjecture stated that every 
link can be
reduced to a trivial link via 3-moves 
(\parbox{2.9cm}{ \psfig{figure=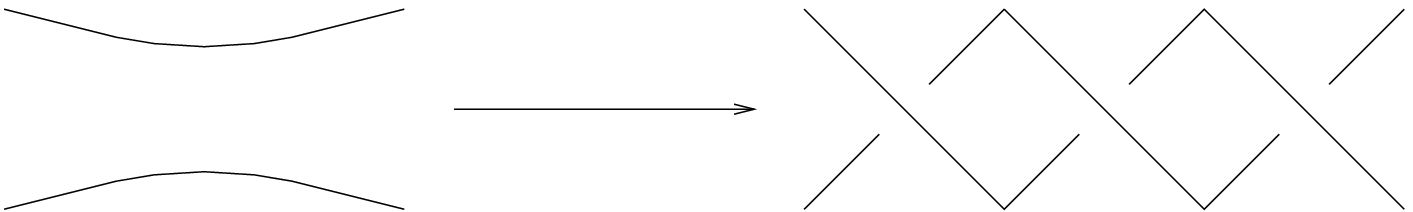,height=0.4cm}}), 
the Nakanishi conjecture stated
that every knot can be reduced to the trivial knot via 4-moves 
(\parbox{3.8cm}{\psfig{figure=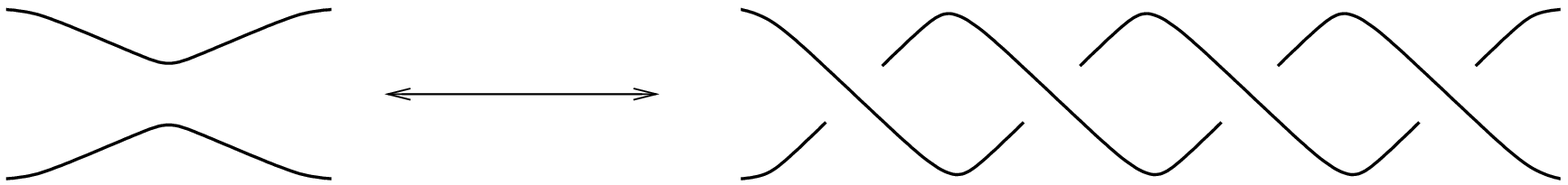,height=0.4cm}}), 
and Harikae-Nakanishi-Uchida conjecture
stated that every link can be reduced to a trivial link via $(2,2)$-moves 
(\parbox{1.6cm}{\psfig{figure=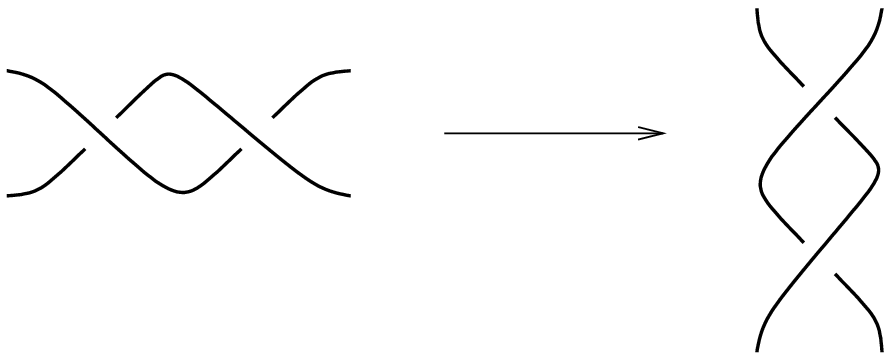,height=0.6cm}}).

It was proven in \cite{D-P-1} that not every link can be reduced to 
a trivial link by 3-moves. Therefore, the Montesinos-Nakanishi conjecture 
does not hold.
The Nakanishi 4-move conjecture is still an open problem. 
In this paper we analyze
 the Harikae-Nakanishi-Uchida conjecture in detail.

\section{$(2,2)$-equivalence of links}


 After 2-,3-, and 4-moves one is tempted to ask 
 about reductions of links by
 5-moves. However,
 it is easy to show that not every link is $5$-move equivalent 
to a trivial link{\footnote{We say that two links are 5-move equivalent
 if one can reach one from the other by a finite number of 5-moves.}}. One
can show, using the Jones polynomial, that the figure eight knot is not
5-move equivalent to a trivial link\footnote{If two unoriented links 
are $5$-move equivalent then their Jones polynomials  for $t^5=-1, t\neq -1$ 
(with any orientations chosen for links) are equal up to an invertible 
element in $Z[t^{\pm 1}]$. 
Furthermore, $V_{T_n}(e^{\pi i/5})=
(-t^{1/2}-t^{-1/2})^{n-1} = (\frac{-1-\sqrt{5}}{2})^{n-1} 
\neq 0$ but $V_{4_1}(e^{\pi i/5})=
t^2-t+1-t^{-1}+t^{-2}= t^{-2}\frac{t^5+1}{t+1}=0$.} 
 (compare \cite{Pr-1}). 
One can, however, introduce a more 
delicate move, called $(2,2)$-move 
(\parbox{1.6cm}{\psfig{figure=dip22move.eps,height=0.6cm}})
such that a $5$- move is a combination of a $(2,2)$-move and its mirror
image, $(-2,-2)$-move 
(\parbox{1.8cm}{\psfig{figure=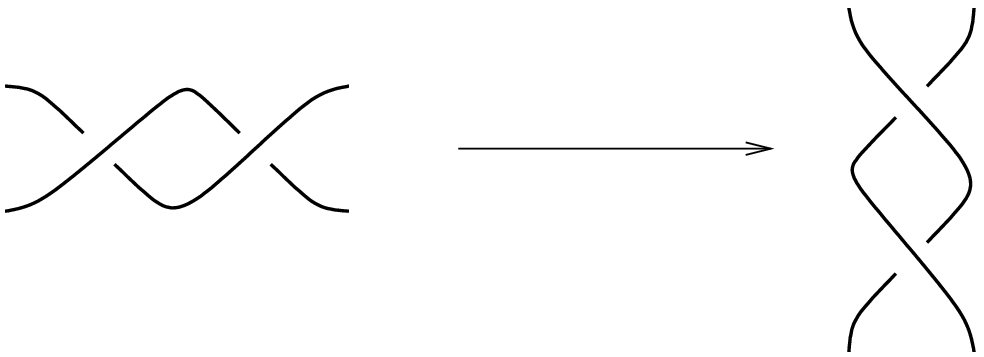,height=0.6cm}}), as
illustrated in Figure 
\ref{fig:dip5-move22-moves} \cite{H-U,Pr-3}. \newline

\begin{figure}
\centerline{
\psfig{figure=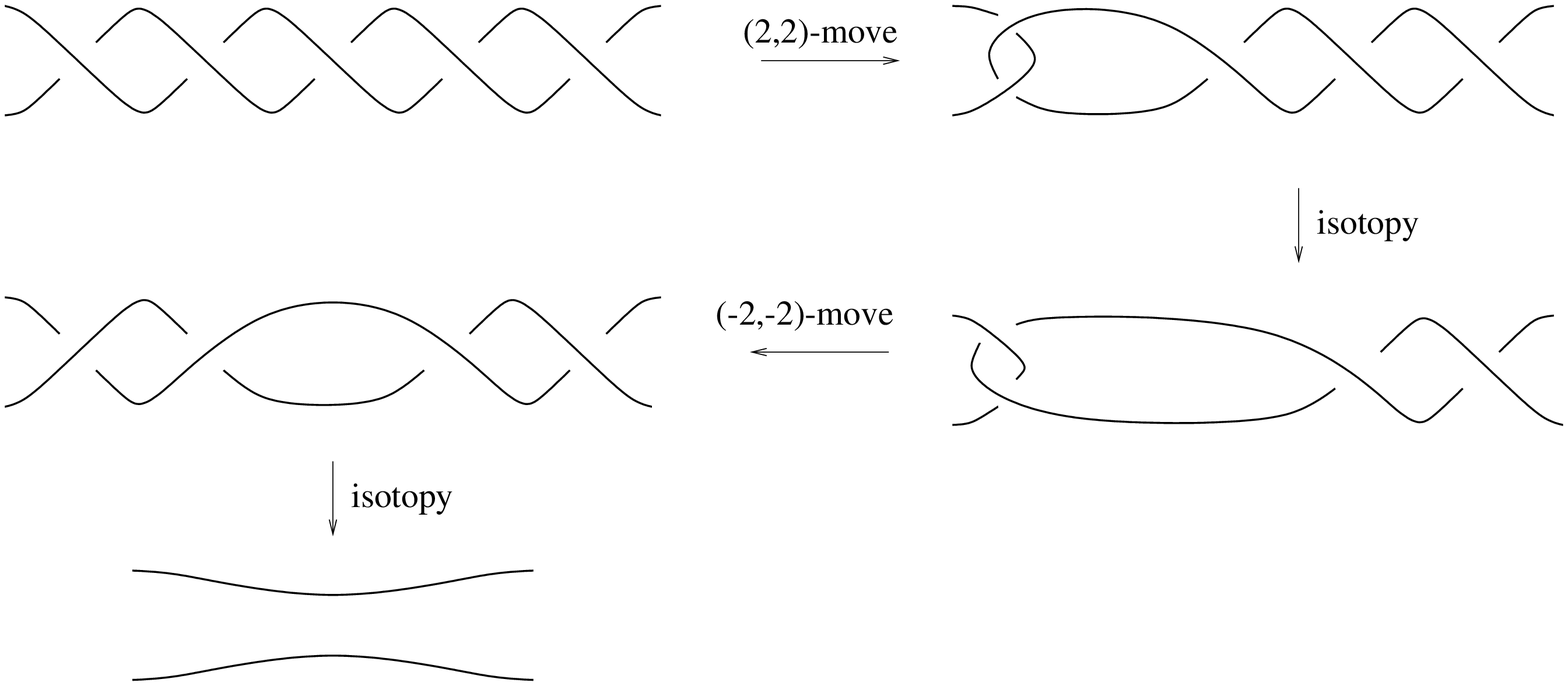,height=5.2cm}
} %
\caption{5-move $\longrightarrow$ (2,2)-moves}\label{fig:dip5-move22-moves}
\end{figure}

We say that two links are $(2,2)$-move equivalent if one can pass 
from one to the other by a finite number of $(2,2)$ and $(-2,-2)$-moves. 
\begin{conjecture}[Harikae, Nakanishi, Uchida 1992 \cite{Kir}]
\label{4.1}\ \newline
Every link is $(2,2)$-move equivalent to a trivial link.
\end{conjecture}

It was shown in \cite{D-P-2} that the conjecture, as stated, does not hold.
The knot $9_{49}$ is a counterexample.
We will show, however, in this paper, to what extent the conjecture holds, 
and how can we modify it so it holds in its full generality.
In particular, we show that every link up to 9 crossings is 
(2,2)-move equivalent to a trivial link, $9_{40}, 9_{49}$ or 
their mirror images. 
The main result of this paper is a classification of closed 
3-braids up to (2,2)-move equivalence.
This result motivates the conjecture that every link can be reduced 
to a trivial link by $\pm(2,2)$-moves and the 
$({\sigma}_1{\sigma}_2)^6$-moves. (See Figure \ref{fig:dipa}.)

\begin{figure}
\centerline{
\psfig{figure=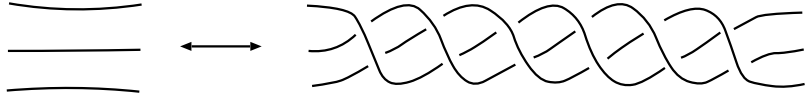,height=1.5cm}
} %
\caption{$({\sigma}_1{\sigma}_2)^6$-move}\label{fig:dipa}
\end{figure}

For 3-braids we prove the following classification theorem.
\begin{theorem} Every closed 3-braid is (2,2)-move equivalent to a 
trivial link or to one of the torus links of the type (3,6), (3,12), (3,18)  
or  (3,24), that is to closures of the 3-braids 
$({\sigma}_1{\sigma}_2)^6$, $({\sigma}_1{\sigma}_2)^{12}$, 
$({\sigma}_1{\sigma}_2)^{18}$  or  $({\sigma}_1{\sigma}_2)^{24}$, respectively.
\end{theorem}

The braid $({\sigma}_1{\sigma}_2)^{30}$ is 5-move equivalent to 
the identity 3-braid (Proposition 2.7.(iii)).

\par\noindent
{\it{Sketch of the proof of the Theorem 2.2.}}
First we show that Harikae-Nakanishi-Uchida conjecture holds for 
algebraic links (Definition 2.3), then we prove the conjecture for 
links up to 8 crossings and eventually, for links up to 9 crossings 
with an exception of (2,2)-equivalence classes of $9_{40}$, 
$9_{49}$, ${\bar{9}}_{40}$  and  ${\bar{9}}_{49}$ (here ${\bar{L}}$ 
denotes the mirror image of $L$). Further, we identify classes of 
$9_{40}$, $9_{49}$ and their mirror images as closures of powers 
of the center of $B_3$. In particular, the knot ${\bar{9}}_{49}$ 
is $(2,2)$-move equivalent to the closure of the 3-braid 
$({\sigma}_1{\sigma}_2)^6$.

We define below $n$-algebraic tangles generalizing a concept of an 
algebraic tangle introduced by J.Conway (his algebraic tangle is 
2-algebraic in our definition)
{\footnote{Here we use only $2$-algebraic tangles and links 
but in \cite{P-Ts1} we proved  that every $3$-algebraic 
link is $3$-move equivalent to a trivial link. 
We hope to address, in the future, the problem of (2,2)-equivalence 
classes of 3-algebraic links.}}.
\\

\begin{definition}[\protect\cite{P-Ts1}]
\label{4.7}
\begin{enumerate}
\item[(i)] $n$-algebraic tangles is the smallest family of n-tangles which
satisfies:
\newline
(0) Any $n$-tangle with $0$ or $1$ crossing is n-algebraic.\newline
(1) If $A$ and $B$ are n-algebraic tangles then $r^{i}(A)\ast r^{j}(B)$ is
n-algebraic, where $r$ denotes the rotation of a tangle by $\frac{2\pi }{2n}$
angle, and * denotes (horizontal) composition of tangles (compare Figure 5).


\item[(ii)] If a link $L$ is obtained from an $n$-algebraic
tangle by closing its endpoints without introducing any new crossings then 
$L$ is called an $n$-algebraic link. 
\end{enumerate}
\end{definition}

\begin{lemma}\label{4.8}
\begin{enumerate}
\item[(i)] Every $2$-algebraic tangle is $(2,2)$-move equivalent to one of
the six 2-tangles that are shown in Figure \ref{fig:dip6basic22dp}  with 
possible additional
trivial components.

\item[(ii)] Every $2$-algebraic link is $(2,2)$-move equivalent to a trivial
link.
\end{enumerate}
\end{lemma}

\begin{proof}
Part (ii) follows immediately from (i) which we prove by induction. Part(i) 
holds for each tangle with no more than one crossing. Thus we have to show  
that if it
holds for tangles in Figure \ref{fig:dip6basic22dp}, say $A$ and $B$, 
then it holds also for $r^{i}(A)\ast r^{j}(B)$. 
All possible compositions involving tangles 
with a crossing are shown in Figure \ref{fig:dipprod-22-reduce}. 
In Figure 5 the reduction of the
``most difficult" case is shown. Other cases are handled in a similar manner.
\end{proof}

\begin{figure}
\centerline{
\psfig{figure=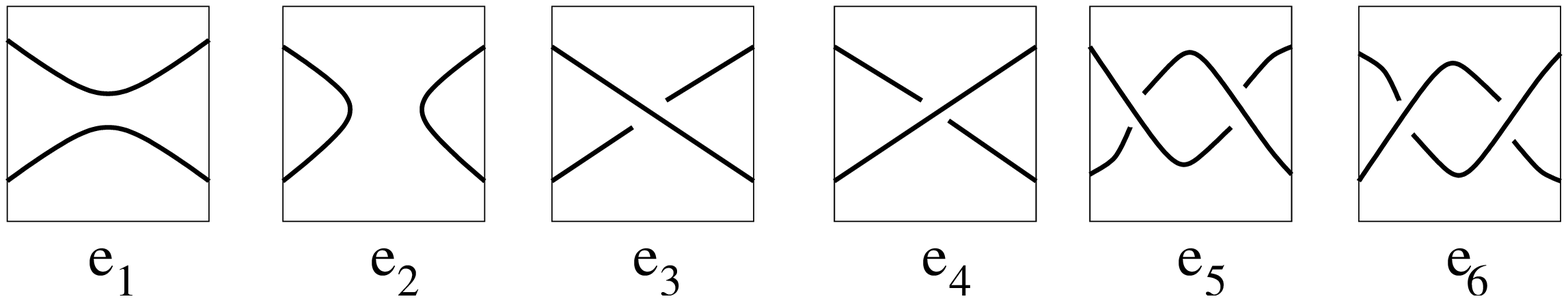,height=2cm} 
} %
\caption{}\label{fig:dip6basic22dp}
\end{figure}

\begin{figure}
\centerline{
\psfig{figure=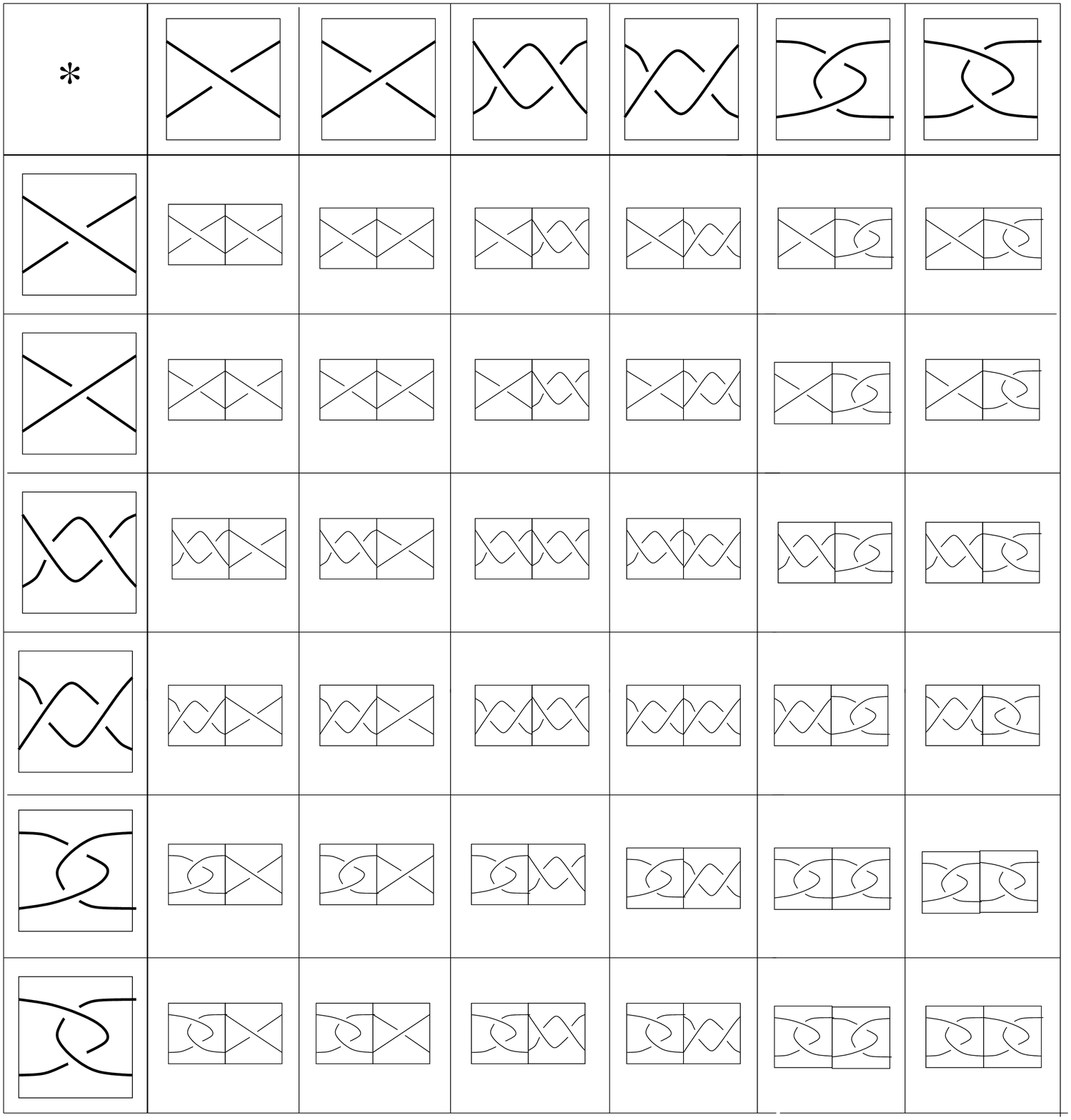,height=15.1cm}} %
\caption{}\label{fig:dipprod-22-reduce}
\end{figure}

\begin{figure}
\centerline{
\psfig{figure=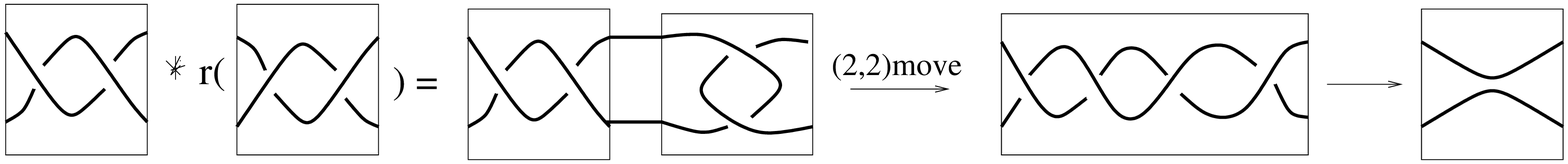,height=1.4cm} 
} %
\caption{}\label{fig:dip22composdp}
\end{figure}

\begin{lemma}\label{4.9}
Every link up to $8$ crossings is $(2,2)$-move equivalent to a 
trivial link.
\end{lemma}
\begin{proof}
According to Conway \cite{Con}, every link up to 8 crossings 
is 2-algebraic with a possible exception of $8_{18}$ 
knot\footnote{To prove
that the knot $8_{18}$ is not 2-algebraic
one considers the 2-fold branched cover of $S^3$ branched along
the knot $M^{(2)}_{8_{18}}$. Montesinos proved that
algebraic knots are covered by Waldhausen graph manifolds \cite{Mo-1}.
Bonahon and Siebenmann showed (\cite{B-S}, Chapter 5)
that $M^{(2)}_{8_{18}}$ is a hyperbolic 3-manifold so it cannot
be a graph manifold.
The knot $9_{49}$ is not 2-algebraic neither because its
2-fold branched cover is a hyperbolic 3-manifold. In
fact, it is the manifold 
I suspected from 1983 to have the smallest volume among
oriented hyperbolic 3-manifolds \cite{I-MPT,Kir,M-V}.
Our work on Burnside groups of links allows us to give a simple 
argument that  the knots $9_{40}$ and $9_{49}$ are not 2-algebraic: 
every algebraic link is (2,2)-move equivalent to a trivial link  
but $9_{40}$ and $9_{49}$ are not. 
However, our method does not work for the knot
$8_{18}$ \cite{D-P-1,D-P-2,D-P-3}.}. The reduction of the $8_{18}$ 
knot to a
trivial link of two components by my students, Jarek Buczy\'{n}ski 
and Mike
Veve, is illustrated in Figure 6.
\end{proof}

\begin{figure}
\psfig{figure=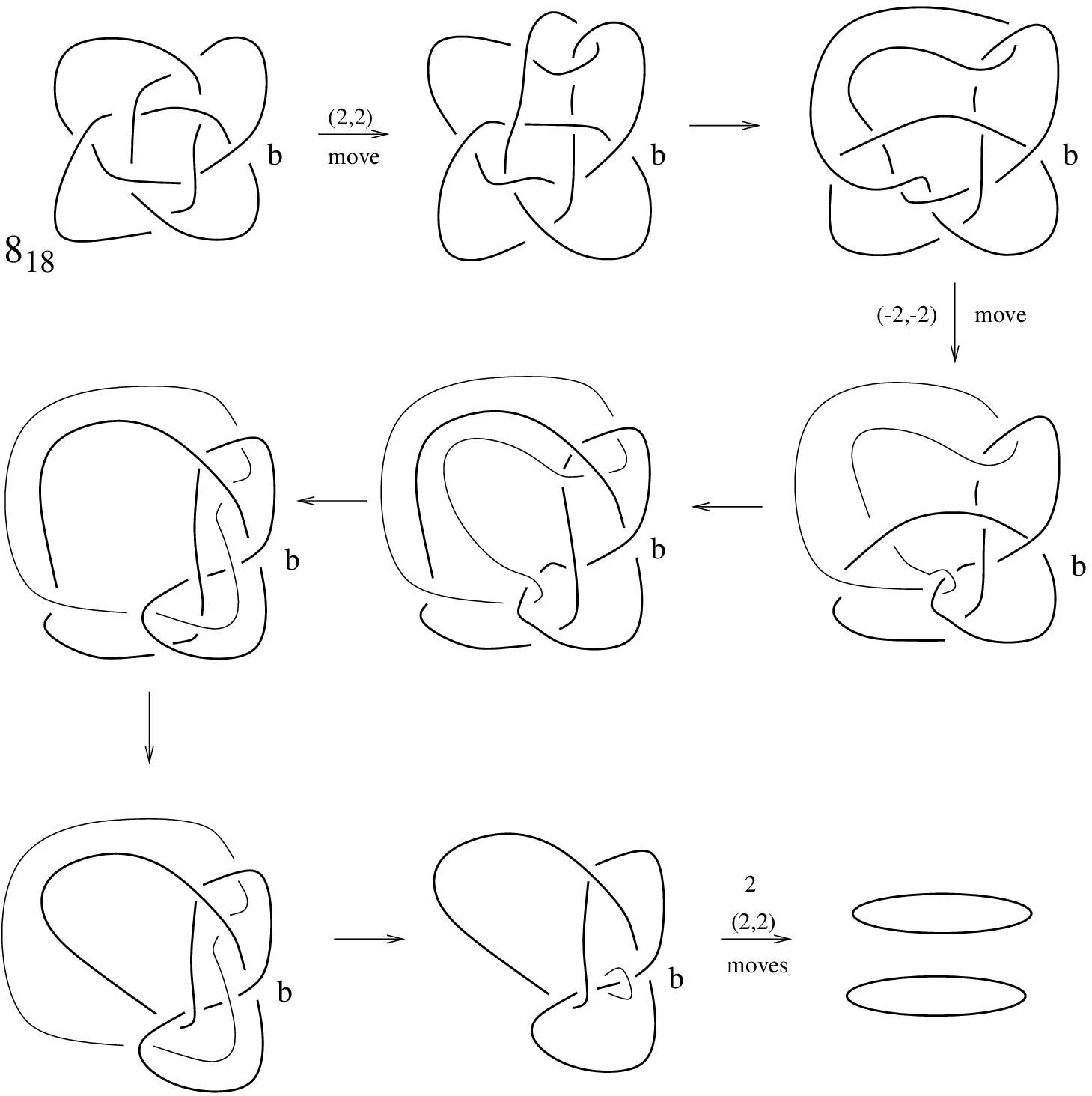,height=13.2cm} 
\caption{Reduction of the $8_{18}$ knot}\label{fig:dip8-18-2-2reduc}
\end{figure}

\begin{lemma}
\label{4.10} Every link up to $9$ crossings is $(2,2)$-move equivalent to a
trivial link or one of the following: 
$9_{40},9_{49},{\bar{9}}_{40},{\bar{9}}_{49}$.
\end{lemma}

\begin{proof}
According to Conway \cite{Con} (checked for us by Slavik Jablan)
the only possible 
$9$ crossing non-algebraic links (up to mirror images) are: $9_{34}$, 
$9_{39}$, $9_{40}$, $9_{41}$, 
$9_{47}$, $9_{49}$, $9_{40}^{2}$,$9_{41}^{2}$,$9_{42}^{2}$,$9_{61}^{2}$. 
We reduce them one by one as it is required in the theorem.

\end{proof}
We proved in \cite{D-P-2} that the knots $9_{40}$ and $9_{49}$ are not 
(2,2)-equivalent to trivial links.

To prove Theorem 2.2, we use the above lemmas and the classical result
 of Coxeter that the group $B_3/({{\sigma}_1}^5)$ is finite  \cite{Cox}. 
More precisely, the group has $600$ elements and there are $45$ conjugacy 
classes. At least 36 of them have  representatives of length at most 8 
(as checked using the computer algebra software GAP). 
By Lemma 2.5, the closures of these braids are $(2,2)-$move
equivalent to trivial links. 
Therefore, it suffices to analyze the remaining nine conjugacy
 classes of $B_{3}/(\sigma _{1}^{5}) $ 
of the length at least $9$. Four of them are listed in the Theorem 2.2:
 ($({\sigma}_1{\sigma}_2)^6$, $({\sigma}_1{\sigma}_2)^{12}$, 
$({\sigma}_1{\sigma}_2)^{18}$, $({\sigma}_1{\sigma}_2)^{24}$). 
The closure of the remaining 5 braids can be reduced to trivial links 
with some effort. 

It follows from \cite{D-P-2} that 
the closures of braids 
$({\sigma}_1{\sigma}_2)^6$, $({\sigma}_1{\sigma}_2)^{12}$, 
$({\sigma}_1{\sigma}_2)^{18}$, $({\sigma}_1{\sigma}_2)^{24}$ 
are not (2,2)-move equivalent to trivial links. 
We do not know, however, whether they are (2,2)-move equivalent among 
themselves. 
The method of \cite{D-P-1,D-P-2}, which uses Burnside groups of links 
does not allow us to separate them (all four links have the same 
5th Burnside group). 

To connect Theorem 2.2. and Lemma 2.6. we prove the following.

\begin{proposition}
\begin{enumerate}
\item[(i)] The knot ${\bar{9}}_{49}$ is $(2,2)$-move equivalent to the
closure of the 3-braid $({\sigma}_1{\sigma}_2)^6$.
\item[(ii)] The knot $9_{40}$ is $(2,2)$-move equivalent to the
closure of the 3-braid $({\sigma}_1{\sigma}_2)^{12}$.
\item[(iii)]  The 3-braid $({\sigma}_1{\sigma}_2)^{30}$ (considered 
as a 3-tangle) is 5-move equivalent to the trivial 3-braid.
\item[(iv)] The knot ${\bar{9}}_{40}$ is $(2,2)$-move equivalent to the
closure of the 3-braid $({\sigma}_1{\sigma}_2)^{18}$.
\item[(v)] The knot $9_{49}$ is $(2,2)$-move equivalent to the
closure of the 3-braid $({\sigma}_1{\sigma}_2)^{24}$.
\end{enumerate}
\end{proposition}

\begin{proof}
\begin{enumerate}
\item[(i)]
The knot $9_{49}$ is $(2,2)$-move equivalent to the
closure of the 3-braid $\alpha_1=(\sigma _{1}^{2}\sigma
_{2}^{-1})^{3}$ (representing $9^2_{40}$)  
(Figure \ref{fig:dip9-49to9-2-40}
).
Furthermore the mirror image of $\alpha_1$ is $5$-move equivalent 
to $(\sigma _{1}\sigma _{2})^{6}$ which is the square of the center of $B_3$. 
We have\\ 
$ \bar\alpha_1 =
(\sigma _{1}^{-2}\sigma _{2})^{3}=(\sigma _{1}^{-2}\sigma _{2})^{2}\sigma
_{2}^{-1}\sigma _{1}^{-1}(\sigma _{2}^{-1}\sigma _{1}\sigma _{2})\sigma
_{1}^{-1}\sigma _{2}=\sigma _{1}^{-2}\sigma _{2}\sigma _{1}^{-2}\sigma
_{1}^{-1}\sigma _{2}^{-1}(\sigma _{1}\sigma _{2}\sigma _{1}^{-1})\sigma
_{2}= $\newline
$\sigma _{1}^{-2}\sigma _{2}\sigma _{1}^{-3}\sigma _{2}^{-2}\sigma
_{1}\sigma _{2}^{2}$.\ 
Further, using three 5-moves we obtain:\\ 
 $\sigma _{1}^{-2}\sigma _{2}\sigma
_{1}^{-3}\sigma _{2}^{-2}\sigma _{1}\sigma _{2}^{2} =_5 \sigma
_{1}^{3}\sigma _{2}\sigma _{1}^{2}\sigma _{2}^{3}\sigma _{1}\sigma
_{2}^{2}=(\sigma _{1}\sigma _{2})^{6}=\Delta _{3}^{4}$,
 which is the square
of the center of the 3-braid group. 
\item[(ii)]
%

\begin{figure}
\centerline{
\psfig{figure=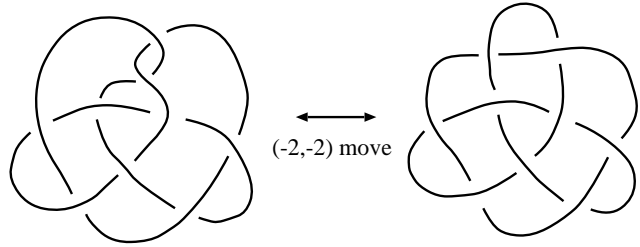,height=3.5cm} 
} %
\caption{$9_{49}$ (2,2)-move reduced to a 3-braid}\label{fig:dip9-49to9-2-40}
\end{figure}

The knot $9_{40}$ is $(2,2)$-move equivalent
to the closure of $\alpha_2=
\sigma_1^2\sigma_2^2\sigma_1^{-2}\sigma_2^2\sigma_1^2\sigma_2^{-2}$ as 
illustrated in Figure \ref{fig:dip9-40-3braid}. 
This is 
conjugated to $(\sigma_1\sigma_2^2\sigma_1^{-2}\sigma_2)
               (\sigma_2\sigma_1^2\sigma_2^{-2}\sigma_1)$. Using a 5-move
twice we get $(\sigma_1\sigma_2^2\sigma_1^{3}\sigma_2)
               (\sigma_2\sigma_1^2\sigma_2^{3}\sigma_1)$. Now we
use expression for a square of the center:
$\Delta_3^4=(\sigma_1\sigma_2)^{6}=
\sigma_1^2\sigma_2^3\sigma_1\sigma_2^2\sigma_1^{3}\sigma_2=
\sigma_2^2\sigma_1^3\sigma_2\sigma_1^2\sigma_2^{3}\sigma_1$
to reduce our braid to $\Delta_3^4\sigma_2^{-3}\sigma_1^{-2}
\sigma_1^{-3}\sigma_2^{-2}\Delta_3^4$ which is reduced further, using two
5-moves, to $\Delta_3^8$ as needed.\\
\item[(iii)]
We have:\ 
$(\sigma_1\sigma_2)^{15}= \sigma_1^3\sigma_2\sigma_1^2\sigma_2\sigma_1^{-1}
(\sigma_1\sigma_2)^{12}=  
 \sigma_1^3\sigma_2(\sigma_1\sigma_2)^3\sigma_1^2(\sigma_1\sigma_2)^3
\sigma_2\sigma_1^{-1}(\sigma_1\sigma_2)^{6}=$  
$$ \sigma_1^3\sigma_2(\sigma_2\sigma_1^2\sigma_2\sigma_1^2)\sigma_1^2
(\sigma_1\sigma_2^2(\sigma_2^2\sigma_1\sigma_2^2\sigma_1)\sigma_1
(\sigma_1\sigma_2^2\sigma_1\sigma_2^2)\sigma_2^2)\sigma_2\sigma_1^{-1}=$$ 
$\sigma_1^3\sigma_2^2\sigma_1^2\sigma_2\sigma_1^5\sigma_2^4\sigma_1
\sigma_2^2\sigma_1^3\sigma_2^2\sigma_1\sigma_2^5\sigma_1^{-1} =_5  $\ 
 $(\sigma_1^{-2}\sigma_2^2)^3 =_{conj}(\sigma_1^{-2}\sigma_2^2)^{-3}=_5
(\sigma_1\sigma_2)^{-15}$ as required.
\item[(iv), (v)]
It follows from (i),(ii), and (iii).
\end{enumerate}
\end{proof}

\begin{figure}
\centerline{
\psfig{figure=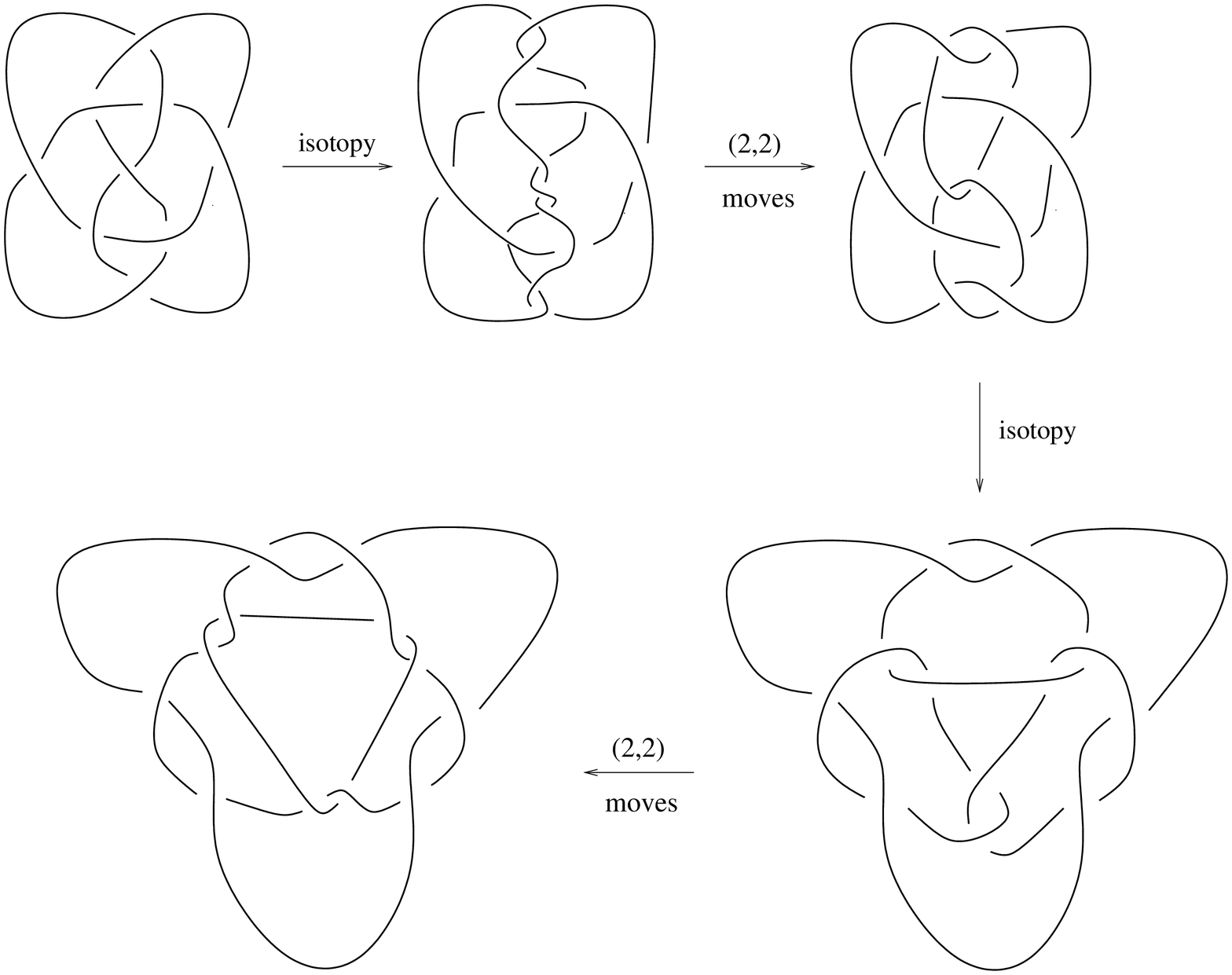,height=9cm} 
} %
\caption{}\label{fig:dip9-40-3braid}
\end{figure}

Let us stress that all four exceptional braids are powers of 
$({\sigma}_1{\sigma}_2)^6$, therefore every closed 3-braids 
can be reduced to a trivial link by $\pm$(2,2)-moves and 
$({\sigma}_1{\sigma}_2)^6$-moves. 
This motivated us to conjecture that every link can be reduced 
to a trivial link by $\pm$(2,2)-moves and  $({\sigma}_1{\sigma}_2)^6$-moves. 
We are glad to announce that the conjecture has been proven recently 
(a week before TWCU Conference) in the joint work with 
T.Tsukamoto  \cite{P-Ts2}.

The $\pm$(2,2)-move is a rational $\pm \frac{5}{2}$-move 
(as illustrated in Figure \ref{fig:dipb2}) and the $n$-move 
is an $\frac{n}{1}$ rational move. We noted in \cite{Pr-4} that the space 
of Fox $n$-colorings ${Col}_n(L)$ (see Definition 3.1) is preserved 
by a rational $\frac{n}{q}$-move for any $q$. 
This can be used to show that different trivial links $U_m$ are 
not (2,2)-move equivalent (because ${Col}_n(U_m)={Z_n}^m$). 
The more sophisticated tool to study rational moves is the non-commutative 
version of Fox $n$-coloring, the $n$th Burnside group of the link $B_n(L)$ 
(see Definition 3.1 (ii)). This group, introduced in \cite{D-P-1,D-P-2}, 
is also invariant under rational $\frac{n}{q}$-moves and has been 
used to prove that knots $9_{40}$ and $9_{49}$ are not 
(2,2)-move equivalent to trivial links.

\section{Kei and (2,2)-moves}

Recall that Kei, \psfig{figure=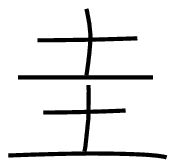,height=0.3cm},
 called also involutive quandle, was introduced by Mituhisa Takasaki 
in 1942 \cite{Tak} as 
an abstract algebra $(Q,*)$ 
with a binary operation $*: Q\times Q \to Q$ satisfying
\begin{enumerate}
\item [(i)] $a*a=a$ for any $a\in Q$  (indempotency condition)
\item [(ii)] $(a*b)*b=a$ (involutory property)
\item [(iii)]$(a*b)*c= (a*c)*(b*c)$ (right distributivity law).
\end{enumerate}
The above axioms have their correspondence in Reidemeister moves 
(see Figure \ref{fig:dipkei-r1r2r3}).\\

\begin{figure}
\centerline{\psfig{figure=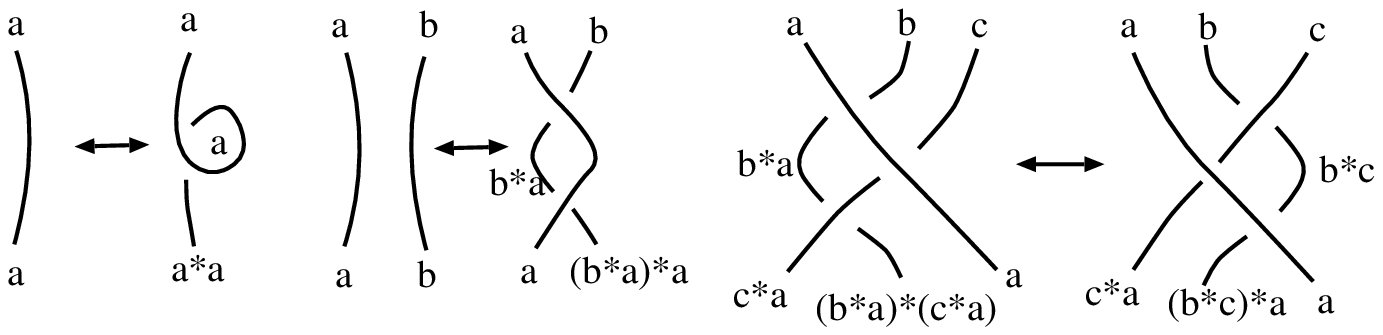,height=3cm}}
\caption{}\label{fig:dipkei-r1r2r3}
\end{figure}

With every group we can associate a core Kei, so that, 
the Kei operation 
is given by $a*b=ba^{-1}b$. 
For every unoriented link diagrams of $L$ we can associate 
the unique 
Kei, ${\bar{Q}}(L)$, by assigning to every arc of the diagram 
a formal variable and taking for every crossing the relation 
as in  Figure \ref{fig:dipd} \cite{Joy}. ${\bar{Q}}(L)$ is a
 link invariant.

\begin{figure}
\centerline{\psfig{figure=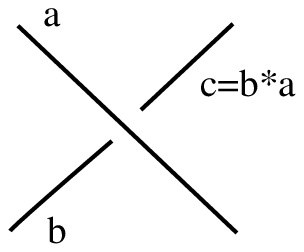,height=2.5cm}}
\caption{}\label{fig:dipd}
\end{figure}

We can distinguish links by comparing their associated Kei. 
For example the Kei of the trivial knot has one element, 
the Kei of the trefoil knot has three elements, the Kei of 
the figure-eight knot has 5 elements and, more generally, 
the Kei of the rational $\frac{n}{q}$ link has  
$n$ elements.

The main examples of Kei that we use are :\\
(i) Dihedral Kei $Z_n$, which is a cyclic group $Z_n$ with the 
operation $i*j \equiv 2j-i \mod n$, and its direct 
sums $Z_n \oplus Z_n \oplus \cdots \oplus Z_n$.
Notice that the operation $i*j \equiv 2j-i$ is an abelian version 
of the product $a*b=ba^{-1}b$.\\
(ii) The free Burnside Kei which is an associated core Kei 
to the Burnside group $B(m,n)$, where $B(m,n)$ is the free 
group $F_m$ divided by relations $w^n=1$ for every $w \in F_m$. \\
Both of these cases can be used to produce invariants of links. 

\begin{definition}
\begin{enumerate}
\item[(i)] 
\begin{enumerate}
\item[(a)] We say that a link (or a tangle) diagram is k-colored if every
arc is colored\footnote{%
We call such a coloring a Fox coloring as R.Fox introduced the construction
when teaching undergraduate students at Haverford College in 1956.} by one
of the numbers $0,1,...,k-1$ $($forming a group $Z_{k})$ in such a way that
at each crossing the sum of the colors of the undercrossings equals 
twice the color of the overcrossing modulo $k$.

\item[(b)] The set of $k$-colorings forms an abelian group, denoted by 
$Col_{k}(D)$. 
\end{enumerate}
\item[(ii)]
\begin{enumerate}
\item[(a)] (\cite{Joy,F-R}) The associated core group of an unoriented link
diagram $D$, $\Pi _{D}^{(2)}$, is the group with generators that correspond
to arcs of the diagram and any crossing $v_{s}$ yields the relation $%
r_{s}=y_{i}y_{j}^{-1}y_{i}y_{k}^{-1}$, where $y_{i}$ corresponds to the
overcrossing and $y_{j},y_{k}$ correspond to the undercrossings at $v_{s}$
(see Figure \ref{fig:dipd}). The core Kei associated to $\Pi _{D}^{(2)}$ is the quotient of the Kei ${\bar{Q}}(D)$ introduced before.

\item[(b)] The unreduced $n$th Burnside group of a link $L$ 
is the quotient
of the associated core group of the the link by its 
normal subgroup 
generated by all relations of the form $w^{n}=1$. 
Succinctly:\ $\hat{B}_{L}(n)=\Pi _{L}^{(2)}/(w^{n})$.

\item[(c)] The $n$th Burnside group of a link is the quotient of the
fundamental group of the double branched cover of 
$S^{3}$ with the link as
the branch set divided by all relations of the form $w^{n}=1$. 
Succinctly:\ $B_{L}(n)=\pi_{1}(M_{L}^{(2)})/(w^{n})$.
\end{enumerate}
\end{enumerate}
\end{definition}

Both groups, the group of Fox $n$-colorings, ${Col}_n(L)$,  
and the Burnside group, $B_L(n)$,  are invariant under 
$n$-moves or, more generally, under rational $\frac{n}{q}$ 
moves \cite{D-P-2}.

This motivates us to analyze general ``behavior" of Kei 
under $n$-moves or, more generally, 
rational moves (including (2,2)-moves). It is useful to notice 
that the free two generator Kei, $Q(2, \infty$), is isomorphic 
to dihedral Kei $Z$ with Kei operation, $ i*j= 2j-i$. 

\begin{lemma} (\cite{Joy})
The Kei homomorphism $\phi : Q(2,\infty) \rightarrow Z$, 
where $\phi$ is given by $\phi(a)=0,\phi(b)=1$, and $a,b$ 
are generators of $Q(2,\infty)$, is an isomorphism. 
\end{lemma}

\par
{\it{Sketch of the Proof:}}
 We have $\phi((w*a)*b)=\phi(w)+2$  and  $\phi(w*a)=-\phi(w)$. 
 These properties imply that $\phi$ is an epimorphism, in particular 
 we have (using the 
left-normed convention: $(...((a*b)*c)*d)...)=a*b*c*d...$) that 
$\phi(b*a)=-1,
\phi(a*b)=2,\phi(a*b*a)=-2, \phi(b*a*b)=3,\phi(b*a*b*a)=-3,
\phi(a*b*a*b)=4,
\phi(a*b*a*b*a)=-4, \phi(b*a*b*a*b)=5$. In general, 
$\phi(b*a*b* \cdots *a*b)=2k+1$ $(2k+1$ letters are used), 
$\phi(a*b* \cdots *a*b)=2k$ $(2k$ letters are used),  
$\phi(a*b* \cdots *a*b*a)=-2k$ $(2k+1$ letters are used), 
and $\phi(b*a \cdots *a*b*a)=-(2k+1)$ $(2k+2$ letters are used).

Monomorphism follows from the fact that every element of $Q(2,\infty)$ 
can be written as a product of elements in the left-normed form 
(we use the identity $x*(y*z)=x*z*y*z$ to achieve this, 
compare \cite{Kam}).

Let us denote by $r_n(a,b)$ the relation given by 
${\phi}^{-1}(0)={\phi}^{-1}(n)$.
For example, $r_5(a,b)$ is given by $a=b*a*b*a*b$. 
It follows from the above that any Kei that is invariant of 
$n$-moves 
must satisfy relations of the form 
${\phi}^{-1}(0)={\phi}^{-1}(n)$. 
In particular, for the 
5-move we get $a=babab$ or equivalently $aba=bab$. 
A similar relation for Kei is required for the 
$\pm(2,2)$-move{\footnote{One can prove more 
generally that 
$\frac{n}{q}$-rational move is preserving Kei satisfying 
universal 
relation of the type ${\phi}^{-1}(0)={\phi}^{-1}(n)$. 
For example, (2,2)-move is equivalent to the 
rational $\frac{5}{2}$-move and the Kei is preserved 
if $a=b*a*b*a*b$ as illustrated in Figure \ref{fig:dipb2}.}.

\begin{figure}
\centerline{
\psfig{figure=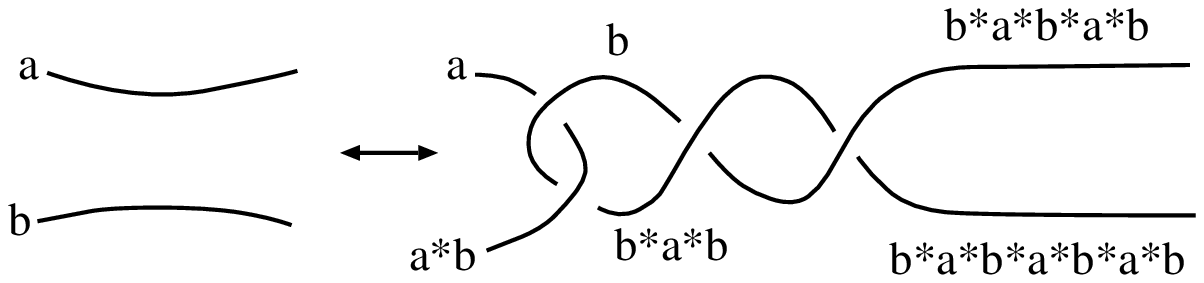,height=2.5cm}}
\caption{}\label{fig:dipb2}
\end{figure}

We denote by $Q(m,n)$ the Kei which has $m$ generators and relations $r_n(u,w)$ for any elements $u,w$. We allow 
$n=\infty$ and then $Q(m, \infty)$ is a free, $m$-generator Kei. 
$Q(m,n)=Q(m,\infty)/(r_n(u,w))$. 
We have a Kei homomorphism $q:Q(m,n) \rightarrow B(m,n)$ being 
the identity on generators. 
To our knowledge, it is an open problem whether $Q(m,3)$ is finite 
or infinite (Takasaki showed that $Q(2,3)$ has 3 elements and 
$Q(3,3)$ has 9 elements. T.Ohtsuki recently (after the TWCU conference) 
computed that $Q(4,3)$ has 81 elements\footnote{The number 
of elements of $Q(4,3)$ also has been computed by 
M.Niebrzydowski. He has also showed that $Q(3,4)$ has $3 \times 2^5=96$ 
elements.}. The result of his calculation is equivalent to showing 
that the Kei homomorphism $Q(4,3) \rightarrow B(4,3)$ is a monomorphism.
Notice that $Q(m,3)$ is a commutative Kei because $r_3=: a=b*a*b$ is 
equivalent to $a*b=b*a$.  
\par
For the analysis of (2,2)-moves, the Kei $Q(m,5)$ is 
very important. 
Every Kei, which produces an invariant of links which 
is preserved 
by (2,2)-moves, is the quotient of $Q(m,5)$. In fact the proof 
that $9_{40}$ and $9_{49}$ are not (2,2)-move equivalent 
to trivial 
links can be formulated in the language of Kei being quotients 
of  $Q(m,5)$. The 5th Burnside group defined in \cite{D-P-2} 
satisfies $a*b*a=b*a*b$. 
We define the $n$th Burnside Kei of a link $L$, ${BQ}_n(L)$, 
as a quotient 
of a link Kei ${\bar{Q}}(L)$ by all relations $r_n(u,w)$ 
(shortly, 
${BQ}_n(L)={\bar{Q}}(L)/(r_n(u,w))$). For the trivial 
link of $m$ 
components $U_m$, we have ${BQ}_n(U_m)=Q(m,n)$. 

\begin{lemma}
${BQ}_n(L)$ is preserved by $n$-moves or, more generally, 
by rational $\frac{n}{q}$-moves.
\end{lemma}

\begin{proof}
The above lemma follows from Lemma 3.2 and the facts 
(essentially known 
to Conway, see \cite{Pr-4}) that the rational 
$\frac{n}{q}$-move 
preserves the space of Fox $n$-colorings, ${Col}_n(L)$.
\end{proof}

On the more philosophical level, Lemma 3.2 can be used 
to show that many 
facts involving rational tangle moves and Fox colorings can 
be formulated in the language of Kei associated to a link.

Motivated by the Burnside's question about finiteness of $B(m,n)$ 
groups we can ask a similar question about Kei $Q(m,n)$. 

\begin{question}\label{3.4}
For which values of $m$ and $n$, is $Q(m,n)$ finite?
\end{question}

$Q(1,n)$ has one element. 
$Q(2,n)$ has $n$ elements and it is isomorphic 
to the dihedral 
Kei $Z_n$. $Q(m,2)$ has $m$ elements and it is a trivial 
Kei \cite{Kam}. 
With support of the knowledge of  
Burnside groups{\footnote{Burnside groups 
of links are instances of groups of finite exponents. 
Our method of analysis of tangle moves rely on the well
developed theory of classical Burnside groups and 
the associated graded Lie rings.
A group $G$ is of a finite exponent if there is a finite integer
$n$ such that $g^{n}=e$ for all $g\in G$. If, in addition, there is no
positive integer $m<n$ such that $g^{m}=e$ for all $g\in G$, then we say
that $G$ has an exponent $n$.
Groups of finite exponents were considered for the first time by
Burnside in 1902 \cite{Bur}. In particular, Burnside himself
was interested in the case when $G$ is
a finitely generated group of a fixed exponent. He asked the question, 
known
as the Burnside Problem, whether there exist infinite and finitely 
generated
groups $G$ of finite exponents. \newline
Let $F_{r}=\langle x_{1},\,x_{2},\,\dots ,\,x_{r}\,|\,-\rangle $ be the 
free
group of rank $r$ and let $B(r,n)=F_{r}/N$, where $N$ is the normal 
subgroup
of $F_{r}$ generated by $\{g^{n}\,|\,g\in F_{r}\}$. 
The group $B(r,n)$ is
known as the $r$th generator Burnside group of exponent $n$. In this
notation, Burnside's question can be rephrased as follows.
For which values of $r$ and $n$ is the Burnside group $B(r,n)$ finite?
$B(1,n)$ is a cyclic group $Z_n$. $B(r,2) = Z_2^n$. Burnside proved 
that $B(r,3)$ is finite for all $r$ and that $B(2,4)$ is finite. 
In 1940 
Sanov proved that $B(r,4)$ is finite for all $r$, and in 1958 M.Hall
proved that $B(r,6)$ is finite for all $r$.
However, it was proved by Novikov and Adjan in 1968
that $B(r,n)$ is infinite whenever $r > 1$, $n$ is odd and $n
\geq 4381$ (this result was later improved by Adjan, who showed
that $B(r,n)$ is infinite if $r > 1$, $n$ is odd and $n \geq 665$).
Sergei Ivanov
proved that for $k \geq 48$ the group $B(2,2^k)$ is infinite.
Lys\"enok found that $B(2,2^k)$ is infinite for $k \geq 13$.
 It is still an open problem, though, whether, 
 for example, $B(2,5)$,
$B(2,7)$ or $B(2,8)$ are infinite or finite \cite{VL,D-P-3}.}
 we could make, maybe too far reaching, prediction 
that $Q(m,n)$ is finite for $n=2, 3, 4$ and $6$.

To show that $9_{40}$ and $9_{49}$ are not (2,2)-equivalent to a trivial 
link we used the $5$-th Burnside group $B_L(5)$ of links. 
Because these groups are Kei's, we could use ${BQ}_L(5)$ to show that fact. 
The Kei ${BQ}_L(5)$ is at least as good as $B_{L}(5)$ in analyzing 
(2,2)-move equivalence. Possibly\footnote{We have been informed by 
M.Niebrzydowski that $BQ_L(5)$ is the same 25 element Kei, 
$Z_5 \oplus Z_5$, for $9_{40}$ and $9_{49}$ knots; 
e-mail; January 14, 2005.}, ${BQ}_L(5)$ can be used to distinguish 
(2,2)-move equivalence classes of $9_{40}, {\bar{9}}_{40},9_{49}$ 
and ${\bar{9}}_{49}$. The first step in this direction is to get 
better understanding of Burnside Kei, ${BQ}_L(5)$. 
In particular to check whether $q:{Q}(5,5) \rightarrow B(5,5)$ has 
a nontrivial kernel.

\section{Application of (2,2)-moves to tangle embedding}\label{3}
We illustrate our method by one example: 
The tangle $T_1$ of Figure \ref{fig:dipc} cannot be embedded 
in the trivial knot,  the Hopf link, and the trefoil knot. 
To demonstrate this, notice that $T_1$ can be reduced by 
two $\pm$(2,2)-moves into the tangle $T_2$ which has a 
trivial component. 
Therefore, every closure of $T_2$ has a nontrivial 
Fox 5-coloring.
The claim in the example follows now from the fact that 
a (2,2)-move is not changing the number of 5-colorings and 
the trivial knot, the Hopf link, and the trefoil knot have 
only trivial Fox 5-colorings.

\begin{figure}
\centerline{
\psfig{figure=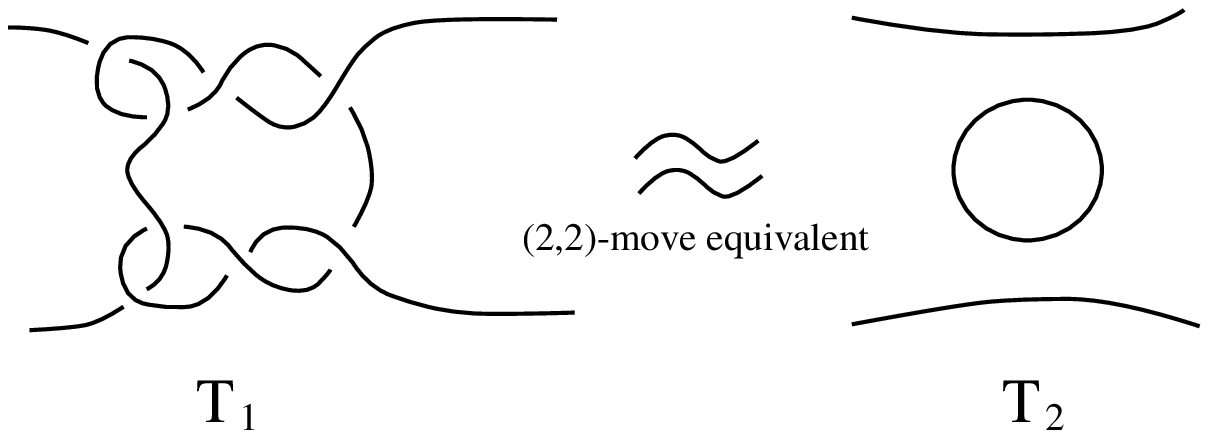,height=3.2cm}}
\caption{}\label{fig:dipc}
\end{figure}

~~
\par
~~
\par
\begin{tabular}{lll}

Mieczys{\l}aw K. D{\c a}bkowski ~~~&Makiko Ishiwata ~~~& J\'ozef H.~Przytycki\\
e-mail: mdab@utdallas.edu ~~~&e-mail: mako@lab.twcu.ac.jp ~~~&e-mail: przytyck@gwu.edu\end{tabular}
\end{document}